\documentclass[11pt]{amsart}

\usepackage{latexsym}
\usepackage{amsmath}
\usepackage{amssymb}
\begin{document}


{\theoremstyle{plain}%
  \newtheorem{theorem}{Theorem}[section]
  \newtheorem{corollary}[theorem]{Corollary}
  \newtheorem{proposition}[theorem]{Proposition}
  \newtheorem{lemma}[theorem]{Lemma}
  \newtheorem{question}[theorem]{Question}
}
{\theoremstyle{remark}
\newtheorem{fact}{Fact}
\newtheorem{remark}[theorem]{Remark}
}

{\theoremstyle{definition}
\newtheorem{definition}[theorem]{Definition}
\newtheorem{example}[theorem]{Example}
}

\newcommand{\bt}{\underline{t}}
\newcommand{\bj}{\underline{j}}
\newcommand{\bi}{\underline{i}}
\newcommand{\dep}{\operatorname{depth}}
\newcommand{\hit}{\operatorname{ht}}
\newcommand{\op}{\overline{P}}
\newcommand{\xab}{X^{\underline{a}}Y^{\underline{b}}}
\newcommand{\Iz}{I_{\Z}}
\newcommand{\Ixp}{I_{\xp}}
\newcommand{\Z}{\mathbb{Z}}
\newcommand{\ax}{\alpha_{\X}}
\newcommand{\bx}{\beta_{\X}}
\newcommand{\kxo}{k[x_0,\ldots,x_n]}
\newcommand{\kx}{k[x_1,\ldots,x_n]}
\newcommand{\popo}{\mathbb{P}^1 \times \mathbb{P}^1}
\newcommand{\pr}{\mathbb{P}}
\newcommand{\pn}{\mathbb{P}^n}
\newcommand{\pnpm}{\mathbb{P}^n \times \mathbb{P}^m}
\newcommand{\pnk}{\mathbb{P}^{n_1} \times \cdots \times \mathbb{P}^{n_k}}
\newcommand{\X}{\mathbb{X}}
\newcommand{\Y}{\mathbb{Y}}
\newcommand{\N}{\mathbb{N}}
\newcommand{\M}{\mathbb{M}}
\newcommand{\Q}{\mathbb{Q}}
\newcommand{\Ix}{I_{\X}}
\newcommand{\pix}{\pi_1(\X)}
\newcommand{\pixt}{\pi_2(\X)}
\newcommand{\pipi}{\pi_1^{-1}(P_i)}
\newcommand{\piqi}{\pi_2^{-1}(Q_i)}
\newcommand{\qpi}{Q_{P_i}}
\newcommand{\pqi}{P_{Q_i}}
\newcommand{\pitk}{\pi_{2,\ldots,k}}
\newcommand{\kdim}{\operatorname{K-}\dim}
\newcommand{\reg}{\operatorname{reg}}
\newcommand{\regI}{\operatorname{reg({\it I}_\X)}}
\newcommand{\ri}{\operatorname{ri}}
\newcommand{\Hx}{\mathcal{H}_{\X}}
\newcommand{\HH}{\mathcal{H}}


\title{The regularity of points in multi-projective spaces}
\thanks{Version: June 9, 2003}
\author{Huy T\`ai H\`a}
\address{Department of Mathematics, University of Missouri, Columbia MO 65211, USA}
\email{tai@math.missouri.edu}
\author{Adam Van Tuyl}
\address{Department of Mathematical Sciences \\
Lakehead University \\
Thunder Bay, ON P7B 5E1, Canada}
\email{avantuyl@sleet.lakeheadu.ca}


\begin{abstract}  
Let $I = \wp_1^{m_1} \cap \ldots \cap \wp_s^{m_s}$ be the defining ideal of a 
scheme of fat points in $\pnk$ with support in generic position. 
When all the $m_i$'s 
are 1, we explicitly calculate the 
Castelnuovo-Mumford regularity of $I$. In general, if at 
least one $m_i \ge 2$, we give an upper bound for the regularity of $I$, 
which extends a result of Catalisano, Trung and Valla.
\end{abstract}

\maketitle


\section*{Introduction}

In this paper, we study the Castelnuovo-Mumford
regularity of defining ideals of sets of points 
(reduced and non-reduced) in a multi-projective space $\pnk$.

If $I \subseteq {\bf k}[x_0,\ldots,x_n]$ is the defining ideal of a projective
 variety $X \subseteq \pr^n$, then the Castelnuovo-Mumford regularity of $I$,
 denoted by $\reg(I)$, is a very important invariant associated to $X$. It 
has been the objective of many authors to estimate $\reg(I)$ since not only 
does it bound
 the degrees of a minimal set of defining equations for $X$, it also 
gives a uniform bound on the degrees of syzygies of $I$. The most fundamental
 situation is when $X$ is a set of points. Examples of work on $\reg(I)$ in
 this case can be seen in \cite{CTV, FL,G, TV}. Recently,
many authors (cf. \cite{CGG,GuMaRa1,Gu,GuVT,VT1}) have been interested
in extending our understanding of points in $\pr^n$ to
sets of points in $\pnk$. We continue this trend
by studying  $\reg(I)$ when $I$ defines a scheme of fat points 
in $\pnk$. 

In the context of $\N^2$-graded rings, Aramova, Crona and De Negri \cite{ACD}
have introduced a finer notion of regularity that places bounds on each
coordinate of the degree of a multi-graded syzygy.   Extending
the definition of regularity to multi-graded rings is also considered in
\cite{HW,MS}.  The usual notion of 
regularity could be treated as a bound on the total degree 
of the multi-graded syzygies. 

The $\N^k$-graded ring
$R = {\bf k}[x_{1,0},\ldots,x_{1,n_1},
\ldots,x_{k,0},\ldots,x_{k,n_k}]$ where $\deg x_{i,j} = e_i$,
the $i^{th}$ basis vector of $\N^k$, is the associated
coordinate ring of $\pnk$.  
Let $\X = \{P_1,\ldots,P_s\}$ be a set of distinct points in $\pnk$.
The defining ideal of $P_i$ is $\wp_i = (L_{1,1},\ldots,L_{1,n_1},\ldots,L_{k,1},\ldots,
L_{k,n_k})$ with $\deg L_{i,j} = e_i$.  If $m_1,\ldots,m_s$
are positive integers, then we want to  study regularity of ideals of the form
$I_Z = \wp_1^{m_1} \cap \cdots \cdots \cap \wp_s^{m_s}$. Such an ideal $I_Z$ defines a scheme of fat points $Z = m_1P_1 + \ldots + m_sP_s$ in $\pnk$. The ideal $I_Z$ is both $\N^k$-homogeneous, and homogeneous
in the normal sense.  Thus, when we refer to $\reg(I_Z)$,
we shall mean its regularity as a homogeneous ideal in $R$,
where $R$ is viewed as an $\N^1$-graded ring.

A set of $s$ points $\X = \{P_1,\ldots,P_s\} \subseteq \pnk$ is said to be 
in generic position if it has maximal Hilbert function 
$H_{\X}(\bi) = \min\{\dim_{\bf k} R_{\bi}, s\}$ for all $\bi \in \N^k$, where 
$R = \bigoplus_{\bi} R_{\bi}$ is the $\N^k$-homogeneous decomposition of $R$.
The existence of such sets is shown in \cite{VT2}.
Our main results consist of explicitly calculating $\reg(I_Z)$ when $Z$ is
in generic position and  reduced 
(i.e. there is no multiplicity at each point), and giving a bound on 
$\reg(I_Z)$ in general. 

In the special case that each $m_i = 1$ and the set
of points is in generic position, we show
\[\reg(I_Z) = \max\{d_1+1,\ldots,d_k+1\}\]
where $d_i = \min\left\{d \in \N ~\left|~ \binom{d+n_i}{d} 
\geq s\right\}\right.$
for each $i = 1,\ldots,k$.  To prove this we use the fact that $I_Z$ is both $\N^k$-homogeneous
and $\N^1$-homogeneous to obtain information about $\reg(I_Z)$.
We also use the Bayer-Stillman criterion for
detecting $m$-regularity \cite{BS}.

We then show that if $\X$ is generic position, and if $m_1 \geq m_2 \geq
\cdots \geq m_s$ with at least one
$m_i \geq 2$, then
\[\reg(I_Z) \leq \max\left\{m_1 + m_2 -1, 
\left[\frac{\sum_{i=1}^s m_i + n_1
-2}{n_1} \right] 
,\ldots,
\left[\frac{\sum_{i=1}^s m_i + n_k
-2}{n_k} \right] \right\} + k. \]
Our strategy is to investigate the regularity index $\ri(R/I_Z)$ of $R/I_Z$, considered as an $\N^1$-graded ring, by
extending the results of \cite{CTV} for fat point schemes in $\pr^n$ to $\pnk$, and then use the fact that $\reg(I_Z) \le \ri(R/I_Z) + k$.

We have organized this papers as follows.  In the first section
we introduce the relevant information about regularity, the regularity
index, and points in multi-projective spaces.  In the second section
we compute the regularity of a defining ideal of a set of points 
in generic position.  In the last section we bound the regularity
for a set of fat points with generic support.

\noindent
{\it Acknowledgments.}
The authors would like to thank A. Conca, for originally poising
this question and some helpful discussions, and E. Guardo,
for her comments on an earlier version of this paper.
This work was begun when the second author visited the Universit\`a di Genova,
and he would like to thank them for their hospitality.  The second
author also acknowledges the financial support of NSERC and INDAM
while working on this project.


\section{Preliminaries}

Throughout this paper ${\bf k}$ denotes an algebraically closed field
of characteristic zero.
In this section, we recall the needed facts about
the Castelnouvo-Mumford regularity, the regularity index,
and points in multi-projective spaces.
Let $S = {\bf k}[x_0,\ldots,x_n]$ be a polynomial ring.

\begin{definition}
A graded $S$-module $M$ is $m$-regular
if there exists a free resolution
\[0 \longrightarrow \bigoplus_j S(-e_{r,j}) \longrightarrow \cdots \longrightarrow
\bigoplus_j S(-e_{1,j}) \longrightarrow \bigoplus_j 
S(-e_{0,j}) \longrightarrow M \longrightarrow 0 \]
of $M$ with $e_{i,j} - i \leq m$ for all $i,j$.  The {\it 
Castelnuovo-Mumford regularity} (or simply, {\it regularity})
of $M$, denoted $\reg(M)$, is the least integer $m$ for which $M$ is
$m$-regular.
\end{definition}

If $I \subseteq S$, then $\reg(I) = \reg(S/I) + 1$.
The  {\it saturation} $\overline{I}$ of the 
ideal $I \subseteq S$ is the ideal
\[\overline{I} := \{F \in S ~|~ 
\mbox{for $i= 1,\ldots,n$, there exists an $r$
such that $x_i^r\cdot F \in I$}\}.\]
$I$ is said to be {\it saturated} if $I = \overline{I}$. The regularity of a saturated ideal does not change if we add a non-zero
divisor.  In fact,

\begin{lemma}[{\cite[Lemma 1.8]{BS}}]  \label{regularnzd}
Let $I \subseteq S$
be a saturated ideal, and suppose $h$ is a non-zero divisor of $S/I$.
Then $I$ is $m$-regular if and only if $(I,h)$ is $m$-regular.
Thus, $\reg(I) = reg((I,h))$.
\end{lemma}
 
The following theorem provides a means to determine
if an ideal is $m$-regular.

\begin{theorem}[{\cite[Theorem 1.10]{BS}} 
Bayer-Stillman criterion for $m$-regularity] \label{bayerstillman}
Let $I \subseteq S$ be an ideal generated in degrees $\leq m$.
The following conditions are equivalent:
\begin{enumerate}
\item[$(i)$] $I$ is $m$-regular.
\item[$(ii)$] There exists $h_1,\ldots,h_j \in S_1$ for some $j \geq 0$
so that 
\begin{enumerate}
\item[$(a)$] $((I,h_1,\ldots,h_{i-1}):h_i)_m = (I,h_1,\ldots,h_{i-1})_m$
for $i = 1,\ldots,j$, and
\item[$(b)$] $(I,h_1,\ldots,h_j)_m = S_m$.
\end{enumerate}
\end{enumerate}
\end{theorem}

The {\it Hilbert function} $H_{M}: \N \longrightarrow \N$
of a graded $S$-module $M$ is defined by
$H_{M}(t) := \dim_{\bf k} M_t$. It is well known (cf. \cite[Theorem 4.1.3]{BH})
that  there exists a unique polynomial $HP_{M}(t)$,
called the {\it Hilbert polynomial of $M$}, such that
$H_{M}(t) = HP_{M}(t)$ for $t \gg 0$.

\begin{definition}
The {\it regularity index} of an $S$-module $M$, denoted $\ri(M)$, is
defined to be
\[\ri(M) := \min\{t ~|~ H_{M}(j) = HP_{M}(j) ~\mbox{for all $j \geq t$}\}.\]
\end{definition}

The regularity and regularity index of an $S$-module are then 
related as follows.

\begin{lemma}[{\cite[Lemma 5.8]{MN}}] \label{regbounds}
If $M$ is a graded $S$-module, then
\[\reg(M) - \dim M + 1 \leq \ri(M) \leq \reg(M) - \operatorname{depth} M + 1.\]
\end{lemma} 
If $M = S/I$, then $\ri(S/I) \leq \reg(S/I) - \operatorname{depth}
S/I + 1 \leq \reg(I)$.  Hence, we have

\begin{corollary} \label{hilbertreg}
If $I \subseteq S$,  then for all $t \geq \reg(I)$,
$H_{S/I}(t) = HP_{S/I}(t)$.
\end{corollary}

Our goal is to investigate $\reg(I)$ when $I$ 
defines either a reduced or non-reduced set of points in $\pnk$
whose support is in generic position.  

Let $R = {\bf k}[x_{1,0},\ldots,x_{1,n_1},\ldots,x_{k,0},\ldots,
x_{k,n_k}]$, with $\deg x_{i,j} = e_i$ where $e_i$ is the $i^{th}$ basis vector of $\N^k$, be the
$\N^k$-graded coordinate ring of $\pnk$. Let $R_{e_i} = {\bf k}[x_{i,0}, \ldots, x_{i, n_i}]$ be the graded coordinate ring of $\pr^{n_i}$ 
for $i = 1,\ldots,k$.
If $P \in \pnk$ is a point, then the ideal $\wp \subseteq R$
associated to $P$ is the prime ideal 
$\wp = (L_{1,1},\ldots,L_{1,n_1},\ldots,L_{k,1},\ldots,L_{k,n_k})$ with $\deg L_{i,j} = e_i$.
Suppose $\X = \{P_1,\ldots,P_s\}$ is a set of distinct points in $\pnk$,
and $m_1,\ldots,m_s$ are $s$ positive integers. Let
\[ I_Z = \wp_1^{m_1} \cap \wp_2^{m_2} \cap \cdots \cap \wp_s^{m_s}\]
where $\wp_i$ is the defining ideal of $P_i$, then $I_Z$ defines a scheme of fat points $Z = m_1P_1 + \ldots + m_sP_s$ in $\pnk$ with support $\X$. When $m_i = 1$ for all $i$, $Z \equiv \X$ is reduced, 
and we usually use $\Ix$ instead of $I_Z$.

Since $\operatorname{ht}(\wp_i) = \sum_{j=1}^k n_j$ for each $i$, it follows that K-$\dim R/{I_Z} = k$.  
Thus, by Lemma \ref{regbounds} we  have 
\[\reg(I_Z) \leq \ri(R/{I_Z}) +k.\]
Note that we have equality if $k = 1$ because then $\operatorname{depth} R/{I_Z} = 1$.
 
We shall find it useful to consider $R/I_Z$ as both 
an $\N^k$-graded ring and as an $\N^1$-graded ring.  
We shall, therefore, use $\HH_{Z}(\underline{t})$ to denote the multi-graded
Hilbert function $\HH_{Z}(\underline{t}) := \dim_{\bf k}
(R/I_Z)_{\underline{t}}$ with $\underline{t} 
= (t_1,\ldots,t_k) \in \N^k$, and $H_{Z}(t)$ to
denote the $\N^1$-graded Hilbert function $H_{Z}:= H_{R/I_Z}$.   
Because $(R/I_Z)_t = 
\bigoplus_{t_1+\cdots+t_k = t} (R/I_Z)_{t_1,\ldots,t_k}$,
we have the  identity:
\[H_{Z}(t) = \sum_{t_1+\cdots + t_k = t} \HH_{Z}(t_1,\ldots,t_k)
\hspace{.25cm}\mbox{for all $t \in \N$.}\]

\begin{definition} A set of $s$ points $\X = \{P_1,\ldots,P_s\}
 \subseteq \pnk$ is said to be in {\it generic position} if 
\[\Hx(\underline{t}) = \min\left\{\dim_{\bf k}R_{\underline{t}}
= \binom{t_1 + n_1}{n_1}\cdots\binom{t_k+n_k}{n_k}, s \right\}
\hspace{.25cm}\mbox{for all $\underline{t} \in \N^k$.}\]
\end{definition}
Further results about points in $\pnk$ can be found in \cite{VT1,VT2}.

\begin{remark}
If $I \subseteq R$ is an $\N^k$-homogeneous ideal, then
the $\N^k$-graded minimal free resolution of $I$ is 
\[ 0 \longrightarrow \mathcal{F}_r \longrightarrow \mathcal{F}_{r-1}
\longrightarrow \cdots \longrightarrow \mathcal{F}_{0} \longrightarrow I \longrightarrow
0\]
where $\mathcal{F}_i = \bigoplus_j 
R(-d_{i,j,1},-d_{i,j,2},\ldots,-d_{i,j,k})$.  Since $I$ is
also homogeneous in the normal sense, the above resolution
also gives a graded minimal free resolution of $I$:
\[ 0 \longrightarrow \mathcal{F}'_r \longrightarrow \mathcal{F}'_{r-1}
\longrightarrow \cdots \longrightarrow \mathcal{F}'_{0} \longrightarrow I \longrightarrow
0\]
where $\mathcal{F}'_i =  \bigoplus_j 
R(-d_{i,j,1}-d_{i,j,2}-\cdots-d_{i,j,k})$ where we view
$R$ as $\N^1$-graded.  So if $I$ is an $\N^k$-homogeneous ideal with
$k \geq 2$, $\reg(I)$ can be interpreted as
a crude invariant that bounds the total
degree of the multi-graded syzygies.
\end{remark}

The following lemma, which generalizes \cite[Lemma 3.3]{VT1}, enables us to 
find non-zero divisors of specific multi-degrees.

\begin{lemma}\label{firstnzd1}
Suppose $\X = \{P_1, \ldots, P_s\}$ is a set of distinct points in $\pnk$, $\wp_1, \ldots, \wp_s$ are the defining ideals of $P_1, \ldots, P_s$, respectively, and $m_1,\ldots,m_s$ are positive integers.  Set $I_Z = \wp_1^{m_1}
\cap \cdots \cap \wp_s^{m_s}$, and fix an $i \in \{1,\ldots,k\}.$  
Then there exists a form $L \in R_{e_i}$ 
such that $\overline{L}$  is a non-zero divisor in $R/{I_Z}$.
\end{lemma}

\section{The regularity of the defining ideal of points in generic position}

Let $\X \subseteq \pnk$ be a set of $s$ reduced points in generic position.
In this section we calculate the Castelnuovo-Mumford regularity
of the defining ideal of $\X$.

For each $i =1,\ldots,k$, set 
$d_i := \min\left\{d ~\left|~ \binom{d+n_i}{d} \geq s 
\right\}\right.$, and let
$D := \max\{d_1 + 1,\ldots,d_k+1\}$.  Note
that if $n_i = \min\{n_1,\ldots,n_k\}$, then $D = d_i + 1$.
Beginning with a combinatorial lemma, we use this notation
to describe the some of the properties of points in generic position.

\begin{lemma}	\label{combinatorics}
Let $n \geq 1$.  Then, for all $a,b \geq 1$,
\[
\binom{a+b+n}{a+b} \leq \binom{a+n}{a}\binom{b+n}{b}.
\]

\end{lemma}

\begin{proof}
Because
\[
\binom{a+b+n}{a+b} 
 =  \frac{(a+b+n) \cdots(a+1+n)}{(a+b)(a+b-1)\cdots(a+1)}\binom{a+n}{a}\]
it is enough to show that the inequality 
\[
\frac{(a+b+n)(a+b-1+n) \cdots(a+1+n)}{(a+b)\cdots(a+1)} \leq \binom{b+n}{b}
\] is true.
This is equivalent to showing that 
\[
\frac{(a+b+n)(a+b-1+n)\cdots(a+1+n)}{(b+n)(b-1+n)\cdots (1+n)}
\leq \frac{(a+b)(a+b-1)\cdots (a+1)}{b(b-1)\cdots 2\cdot1}.
\]
Rewriting the above expression, we see that we need to show that
\[
\left[1 + \frac{a}{b+n}\right]\left[1 + \frac{a}{b-1+n}\right]\cdots
\left[1 + \frac{a}{1+n}\right] \leq
\left[1 + \frac{a}{b}\right]\left[1 + \frac{a}{b-1}\right]\cdots
\left[1 + \frac{a}{1}\right]
.\]
But since 
$\left[1 + \frac{a}{b+n-j}\right] \leq \left[1 + \frac{a}{b-j}\right]$ 
for $j = 0,\ldots,b-1$ we are finished.
\end{proof}

\begin{corollary} 	\label{combinatoricscor}
Let $\X \subseteq \pnk$ be $s$  points in generic position.  
If $(t_1,\ldots,t_k)
\in \N^k$ is such that $t_1 + \cdots + t_k = D - 1$, then
$\Hx(t_1,\ldots,t_k) = s$.
\end{corollary}

\begin{proof} Suppose that $n_i = \min\{n_1,\ldots,n_k\}$,
and hence, $D-1 = d_i$.  Lemma \ref{combinatorics} then gives 
\begin{eqnarray*}
\binom{t_1+n_1}{t_1}\binom{t_2 +n_2}{t_2}\cdots
\binom{t_k+n_k}{t_k} 
& \geq &  
\binom{t_1+n_i}{t_1}\binom{t_2+ n_i}{t_2} \cdots \binom{t_k+n_i}{t_k} 
 \geq 
\binom{d_i + n_i}{d_i}
\end{eqnarray*}
Since $\binom{d_i + n_i}{d_i} \geq s$, we have $\Hx(t_1,\ldots,t_k) = s$.
\end{proof}

Recall that if $m \in \N$, then 
$\binom{t + m}{m}$  denotes
the polynomial
\[\binom{t + m}{m} = \frac{(t+m)(t+(m-1))\cdots (t+1)}{m!}.\]

\begin{proposition}  \label{genprop}
Let $\Ix$ be the defining ideal of $s$ points $\X \subseteq
\pnk$ in generic position.
\begin{enumerate}
\item[$(i)$] As an $\N^1$-graded ideal, $\Ix$ is generated by
forms of degree $\leq D$.
\item[$(ii)$] As an $\N^1$-graded ring, $R/\Ix$ has Hilbert polynomial 
$HP_{R/\Ix}(t) = s\binom{t+k-1}{k-1}.$
\item[$(iii)$] Fix an $i \in \{1,\ldots,k\}$ and let $L$ be
the non-zero divisor of Lemma \ref{firstnzd1} of degree $e_i$.
If $\bt = (t_1,\ldots,t_k) \in \N^k$ is such that
$t_1 + \ldots + t_k \geq D$ and $t_i > 0$, then
$(\Ix,L)_{\bt} = R_{\bt}.$
\end{enumerate}
\end{proposition}

\begin{proof}
For $(i)$ it suffices to show that for all $\bt = (t_1,\ldots,t_k) \in \N^k$
with $t_1 + \cdots + t_k \geq D+1$, $(\Ix)_{\bt}$ contains
no new minimal generators.  If $\bt \in \N^k$ is such a tuple, then
there exists $l,m \in \{1,\ldots,k\}$, not necessarily distinct,
such that $\bt - e_l - e_m \in \N^k$.  By Corollary ~\ref{combinatoricscor}
it follows that $\Hx(\bt-e_l-e_m) = \Hx(\bt-e_l) = s$
since $t_1 + \cdots + t_k -2 \geq D-1$.  Now apply the results of 
\cite{VT2} to conclude that $(\Ix)_{\bt}$ contains
no minimal generators.

Since $\X$
is in generic position, for $t \gg 0$ we have
\[
H_{\X}(t) = \sum_{t_1 + \cdots + t_k = t} \Hx(t_1,\ldots,t_k) =
\sum_{t_1 + \cdots +t_k = t} s = s\binom{t+k-1}{k-1}.
\]
Since $HP_{R/\Ix}$ is the unique polynomial that agrees
with $H_{\X}$ for $t \gg 0$, $(ii)$ now follows.

To prove $(iii)$ we only consider the case $i = 1$.  
Since $\overline{L}$ is a non-zero divisor,
the exact sequence 
\[
0 \longrightarrow (R/\Ix)(-e_1) \stackrel{\times \overline{L}}
{\longrightarrow} R/\Ix \longrightarrow R/(\Ix,L) \longrightarrow 0
\]
implies that
\[
\mathcal{H}_{R/(\Ix,L)}(t_1,\ldots,t_k)
= \Hx(t_1,\ldots,t_k) - \Hx(t_1-1,t_2,\ldots,t_k)
\hspace{.25cm}\mbox{for all $\underline{t} \in \N^k$}
\]
where $\Hx(t_1-1,t_2,\ldots,t_k) = 0$ if $t_1 -1 <0$.
Now suppose that $t_1 + \ldots + t_k \geq D$ with $t_1 > 0$.
Since $(t_1 -1)+ t_2+\cdots+t_k \geq D-1$, by Corollary \ref{combinatoricscor}
we have $\Hx(t_1,\ldots,t_k) = \Hx(t_1-1,t_2,\ldots,t_k) = s$.
Thus $\mathcal{H}_{R/(\Ix,L)}(t_1,\ldots,t_k) = 0$, or
equivalently, $(\Ix,L)_{t_1,\ldots,t_k} = R_{t_1,\ldots,t_k}$. 
\end{proof}

\begin{theorem}
Let $\Ix$ be the defining ideal of $s$ points 
$\X \subseteq \pnk$ in generic position.
Then
\[
\reg(\Ix) = \max\{d_1+1,\ldots,d_k + 1\}
\]
where $d_i := \min\left\{d ~\left|~ \binom{d +n_i}{d} \geq s\right\}\right.$
for $i = 1,\ldots,k$.
\end{theorem}

\begin{proof}
Without loss of generality, 
we  assume that $n_1 \ge n_2 \ge \ldots \ge n_k \ge 1$.  
It thus suffices to show that $\reg(\Ix) = d_k + 1 = 
\max\{d_1+1,\ldots,d_k+1\}$.  

We first show that $\reg(\Ix) > d_k$.
By Lemma \ref{firstnzd1} there is a
non-zero divisor 
$\overline{L}$ of 
$R/\Ix$ with $\deg L = e_k$.
As an $\N^1$-homogeneous element of $R$, $\deg L = 1$.
Since $\Ix$ is saturated, by Lemma \ref{regularnzd} is it
is enough to show $\reg(\Ix,L) > d_k$.

From the  short exact sequence
\[0 \longrightarrow (R/\Ix)(-1) \stackrel{\times \overline{L}}{\longrightarrow}
R/\Ix \longrightarrow R/(\Ix,L) \longrightarrow 0.\]
of $\N^1$-graded rings, and from Proposition \ref{genprop} $(ii)$ 
we deduce that
\[
HP_{R/(\Ix,L)}(t)  = HP_{R/\Ix}(t) - HP_{R/\Ix}(t-1) 
 =  s\binom{t+(k-2)}{k-2}.\]

If we can show that $HP_{R/(\Ix,L)}(d_k) \neq H_{R/(\Ix,L)}(d_k)$,
then by Corollary ~\ref{hilbertreg}, we can conclude that
$\reg(\Ix,L) > d_k$.  So, write $H_{R/(\Ix,L)}(d_k) = A + B$ where
\[
A :=  \sum_{t_1 + \cdots + t_{k-1} = d_k} \mathcal{H}_{R/(\Ix,L)}(t_1,\ldots,t_{k-1}, 0)
~~\mbox{and}~~ 
B  :=  \sum_{t_1 + \cdots + t_k = d_k, ~t_k > 0} \mathcal{H}_{R/(\Ix,L)}(t_1,t_2,\ldots,t_k). 
\]
From the  short exact sequence
\[0 \longrightarrow (R/\Ix)(-e_k) \stackrel{\times \overline{L}}{\longrightarrow}
R/\Ix \longrightarrow R/(\Ix,L) \longrightarrow 0\]
of $\N^k$-graded rings,
we have 
\[\mathcal{H}_{R/(\Ix,L)}(t_1,\ldots,t_k) = 
\mathcal{H}_{R/\Ix}(t_1,\ldots,t_k) - 
\mathcal{H}_{R/\Ix}(t_1,\ldots,t_{k-1}, t_k - 1)\]
where $\mathcal{H}_{R/\Ix}(t_1,\ldots,t_{k-1}, t_k -1) = 0$ if
$t_k = 0$.  Thus 
\[A = \sum_{t_1+\cdots+t_{k-1} = d_k} \mathcal{H}_{R/\Ix}(t_1,\ldots,t_{k-1},0).\]
Since $t_1+\cdots+t_{k-1} = d_k$,
by Corollary ~\ref{combinatoricscor}
we have $\mathcal{H}_{R/\Ix}(t_1,\ldots,t_{k-1}, 0) = s$.
Hence,
\[
A = \sum_{t_1+\cdots +t_{k-1} = d_k} s = s\binom{d_1+k-2}{k-2} =
HP_{R/(\Ix,L)}(d_k).
\]
On the other hand, because 
$d_k = \min\left\{d ~\left|~ \binom{d + n_k}{d} \geq s \right\}\right.$
\begin{eqnarray*}
B &\geq & \HH_{R/(\Ix,L)}(0, \ldots, 0, d_k) = \Hx(0, \ldots, 0, d_k) - 
\Hx(0, \ldots, 0, d_k-1) \\
&=&  s - \binom{d_k-1+n_k}{d_k-1} > 0.\\
\end{eqnarray*}
Thus, $H_{R/(\Ix,L)}(d_k) = HP_{R/(\Ix,L)}(d_k) + B > 
HP_{R/(\Ix,L)}(d_k)$, as desired.

We now show that $\reg(\Ix) \leq d_k + 1$ by demonstrating
that $\Ix$ is $(d_k + 1)$-regular.
By Proposition \ref{genprop} $(i)$, as an $\N^1$-graded
ideal $\Ix$ is generated by elements of degree $\leq d_k + 1$.
For each $i \in \{1,\ldots,k\}$, by Lemma \ref{firstnzd1} there
exists a non-zero divisor $\overline{L}_i \in R/\Ix$ with
$\deg L_i = e_i$. 
After a change of variables in the $x_{1,j}$'s, a change
of variables in the $x_{2,j}$'s, etc.,
we can assume that $L_i = x_{i,0}$ for $i =1,\ldots,k$.

By the Bayer-Stillman criterion (Theorem ~\ref{bayerstillman}),
to show that $\Ix$ is  $(d_k + 1)$-regular, it is enough to prove:
\begin{enumerate}
\item[$(a)$] $((\Ix,x_{1,0},\ldots,x_{j-1,0}):x_{j,0})_{d_k+1} = 
(\Ix,x_{1,0},\ldots,x_{j-1,0})_{d_k+1}$
for $j = 1,\ldots,k$, 
\item[$(b)$] $(\Ix,x_{1,0},\ldots,x_{k,0})_{d_k+1} = R_{d_k+1}$.
\end{enumerate}

\noindent
{\it Proof of $(a)$.}
We need to only show the 
non-trivial inclusion $[(\Ix,x_{1,0},\ldots,x_{j-1,0}):x_{j,0}]_{d_k+1}
\subseteq (\Ix,x_{1,0},\ldots,x_{j-1,0})_{d_k+1}$ for each $j$.
If $j = 1$, then the statement holds because $x_{1,0}$ is a non-zero
divisor.

So, suppose $j > 1$.   Set $J :=  [(\Ix,x_{1,0},\ldots,x_{j-1,0}):x_{j,0}]$.
Because $J$ is also $\N^k$-homogeneous, if $F \in J_{d_k+1}$,
then can assume that $\deg F = \bt = (t_1,\ldots,t_k)$
with $t_1 + \cdots + t_k = d_k + 1$.  There are now two cases to consider.

In the first case, one of $t_1,\ldots,t_{j-1} > 0$.
Suppose $t_l >0$ with $1 \leq l \leq (j-1)$.  
Then by Proposition \ref{genprop} $(iii)$ we have
$F \in R_{\bt} \subseteq (\Ix,x_{l,0})_{\bt}
\subseteq (\Ix,x_{1,0},\ldots,x_{j-1,0})_{\bt}.$
Since 
$(\Ix,x_{1,0},\ldots,x_{j-1,0})_{\bt} \subseteq
(\Ix,x_{1,0},\ldots,x_{j-1,0})_{d_k+1}$ (as vector spaces), we are finished.

In the second case, $t_1 = t_2 = \cdots = t_{j-1} = 0$.
Then $Fx_{j,0} \in 
(\Ix,x_{1,0},\ldots,x_{j-1,0})_{0,\ldots,0,t_j+1,\ldots,t_k}$.
But since 
\[(\Ix,x_{1,0},\ldots,x_{j-1,0})_{0,\ldots,0,t_j+1,\ldots,t_k}
= (\Ix)_{0,\ldots,0,t_j+1,\ldots,t_k},\] we have $Fx_{j,0}
\in (\Ix)_{0,\ldots,0,t_j+1,\ldots,t_k}$.  But because 
$x_{j,0}$ is a non-zero divisor of $R/\Ix$, 
\[F \in (\Ix)_{0,\ldots,0,t_j,\ldots,t_k}
\subseteq (\Ix,x_{1,0},\ldots,x_{j-1,0})_{0,\ldots,0,t_j,\ldots,t_k}
\subseteq (\Ix,x_{1,0},\ldots,x_{j-1,0})_{d_k+1}.\]

\noindent
{\it Proof of $(b)$.}  
Since 
$
R_{d_k+1} = \bigoplus_{t_1+\cdots+t_k = d_k +1} R_{t_1,\ldots,t_k}$
and because $(\Ix,x_{1,0},\ldots,x_{k,0})$ is also $\N^k$-homogeneous,
it is enough to show that 
$R_{\bt} \subseteq (\Ix,x_{1,0},\ldots,x_{k,0})_{\bt}$
for all $\bt = (t_1,\ldots,t_k) \in \N^k$ 
with $t_1 + \cdots + t_k = d_k + 1$.
But for any $\bt \in \N^k$ with $t_1 + \cdots + t_k = d_k + 1$,
there exists at least one $t_l > 0$. Thus, by
Proposition \ref{genprop} $(iii)$ we have $R_{\bt} \subseteq 
(\Ix,x_{l,0})_{\bt} \subseteq 
(\Ix,x_{1,0},\ldots,x_{k,0})_{\bt}$, thus completing the proof of $(b)$.

Since we have just shown $d_k < \reg(\Ix) \leq d_k + 1$,
the desired conclusion now follows.
\end{proof}

\begin{remark}
If $\X$ is a set of $s$ points in generic position in $\pr^n$,
we recover the well known result that $\reg(\Ix) = d + 1$
where $d = \min\{l ~|~ \binom{l+n}{n} \geq s\}$.
\end{remark}


\section{Bounding the regularity of fat points in $\pnk$}

Let $\X = \{P_1,\ldots,P_s\} \subseteq \pnk$ and 
$m_1 \geq \cdots \geq m_s \in \N^+$.  Suppose $\wp_i$ is 
the defining ideal of $P_i$ for $i = 1, \ldots, s$. Let $I = I_Z = \wp_1^{m_1} \cap \cdots \cap \wp_s^{m_s}$. In this section, we give an upper bound for $\reg(I)$ when $\X$ is in generic position. If we consider $R/I$ as
an $\N^1$-graded ring, then by Lemma \ref{regbounds} 
\[\reg(I) = \reg(R/I)  + 1 
\leq \ri(R/I) + \dim R/I = \ri(R/I) + k.\]
To bound $\reg(I)$, it is therefore enough to
bound  $\ri(R/I)$. For convenience, we assume that $n_1 \ge \ldots \ge n_k$.
In the sequel, we shall also abuse notation by writing
$L$ for the form
$L \in {\bf k}[x_{j,0},\ldots,x_{j,n_j}]$, the hyperplane
$L$ in $\pr^{n_j}$ defined by $L$, and the subvariety of $\pnk$ 
defined by $L$.

\begin{lemma} \label{riforpoint}
If $\wp$ is the defining ideal of  point $P \in \pnk$, then
\[\ri(R/\wp^a) = a - k  ~~\mbox{for all $a \geq 1$}.\]
\end{lemma}

\begin{proof} Since $\wp$ defines a complete intersection
of  height $\sum_{i=1}^k n_i$,  Lemma \ref{regbounds}
gives $\ri(R/\wp^a) = \reg(R/\wp^a)-k+1$.  The conclusion
follows since $\reg(\wp^a) = a\reg(\wp) = a$ by 
\cite[Theorem 3.1]{CH}.
\end{proof}

\begin{lemma} \label{riformanypoints}
Suppose $P_1,\ldots,P_r,P$ are points in generic position in $\pnk$,
and let $\wp_i$ be the defining ideal of $P_i$ and let $\wp$
be the defining ideal of $P$.  Let $m_1,\ldots,m_r$ and $a$
be positive integers, $J = \wp_1^{m_1} \cap \cdots \cap \wp_r^{m_r}$,
and $I = J \cap \wp^a$.  Then
\[\ri(R/I) \leq \max\left\{ a-k, \ri(R/J), \ri(R/(J+\wp^a))\right\}.\]
Furthermore, $R/(J+\wp^a)$ is artinian.
\end{lemma}

\begin{proof}
The short exact sequence of $\N^1$-graded rings
\[0 \longrightarrow R/I \longrightarrow
R/J \oplus R/\wp^a \longrightarrow R/(J+\wp^a) \longrightarrow 0\]
yields $H_{R/I}(t) = H_{R/J}(t) + H_{R/\wp^a}(t) -
H_{R/(J+\wp^a)}(t)$.  Combining this with Lemma \ref{riforpoint} gives
\[\ri(R/I) \leq \max\left\{a-k,\ri(R/J), \ri(R/(J+\wp^a))\right\}.\]

To show that $R/(J+\wp^a)$ is artinian, we need to show that
there exists $b$ such that for all $\bt = (t_1,\ldots,t_k) \in \N^k$,
if there is $t_j \geq b$, then  $(R/(J+\wp^a))_{\bt} = 0$.
So, it suffices to show that there exists such a $b$ so
that for all $\bt = (t_1,\ldots,t_k)$ with $t_j \geq b$ for some
$j$, then all monomials of $R$ of degree $\bt$
are in $(J+\wp^a)$.  Suppose $M$ is a monomial
in $R$ of degree $\bt$.  Then $M = N_1N_2\cdots N_k$ where
$N_l$ are monomials in $\{x_{l,0},\ldots,x_{l,n_l}\}$
and of degree $t_l$.  It is enough to show $N_j \in (J+\wp^a)$.

Let $Q_1,\ldots,Q_r,Q$ be the projections of $P_1,\ldots,P_r,P$ in $\pr^{n_j}$.
Since the points are in generic position, the projections are distinct.
Let $\mathcal{Q}_1,\ldots,\mathcal{Q}_r$ and $\mathcal{Q}$ be the defining
ideals of $Q_1,\ldots,Q_r,Q$ in $A = {\bf k}[x_{j,0},\ldots,x_{j,n_j}]$.
Then it is easy to see that 
$A/(\mathcal{Q}_1^{m_1} \cap \cdots \cap \mathcal{Q}_r^{m_r} + \mathcal{Q}^a)$
is artinian.  As well, $\mathcal{Q}_1^{m_1} \cap \cdots \cap 
\mathcal{Q}_r^{m_r} \subseteq J$ and $\mathcal{Q}^a \subseteq \wp^a$,
and thus $\mathcal{Q}_1^{m_1} \cap \cdots \cap \mathcal{Q}_r^{m_r} + 
\mathcal{Q}^a \subseteq (J+\wp^a)$, and this is what needs to be
shown.
\end{proof}

From Lemma \ref{riformanypoints}, to estimate $\ri(R/I)$
we need to estimate $\ri(R/(J+\wp^a))$, or equivalently,
the least integer $t$ such that $(R/(J+\wp^a))_t = 0$,
when this ring is consider as $\N^1$-graded.

\begin{lemma} \label{bounds}
With the same hypotheses as in Lemma \ref{riformanypoints},
and considering the $\N^1$-gradation, we have
\begin{enumerate}
\item[$(i)$] $H_{R/(J+\wp^a)}(t) = \sum_{i=0}^{a-1} \dim_{\bf k}
\left[(J+\wp^i)/(J+\wp^{i+1})\right]_t ~~\mbox{for all $t \geq 0$}$.
\item[$(ii)$] If $P = [1:0:\cdots :0] \times \cdots \times [1:0:\cdots
:0]$ then $\left[(J+\wp^i)/(J+\wp^{i+1})\right]_t = 0$ if and only if
either $i > t$, or $i < t$ and $GM \in (J+ \wp^{i+1})$ for every monomial
$M$ of degree $i$ in $\{x_{1,1},\ldots,x_{1,n_1},\ldots,x_{k,1},\ldots,
x_{k,n_km}\}$, and every monomial $G$ of degree $t -i $
in $\{x_{1,0},x_{2,0},\ldots,x_{k,0}\}$.
\end{enumerate}
\end{lemma}

\begin{proof}
The first assertion follows from the short exact sequences:
\[ 0 \longrightarrow (J + \wp^i)/(J+ \wp^{i+1}) \longrightarrow
R/(J + \wp^{i+1}) \longrightarrow R/(J + \wp^{i})
\longrightarrow 0 \]
where $i = 0, \ldots, a-1$.

To prove $(ii)$, if $i > t$, then 
$(J + \wp^i)_t = (J + \wp^{i+1})_t = J_t$.  So suppose $i < t$. 
We see that $\wp = (x_{1,1},\ldots,x_{1,n_1},
\ldots,x_{k,1},\ldots,x_{k,n_k})$.  Thus
$((J+\wp^i)/(J+\wp^{i+1}))_t = 0$ if and only if $(\wp^i)_t 
\subseteq (J + \wp^{i+1})_t$ if and only if $FM \in 
(J + \wp^{i+1})$ for every monomial $M$ of degree 
$i$ in $\{x_{1,1},\ldots,x_{1,n_1},\ldots,x_{k,1},\ldots,x_{k,n_k}\}$
and every form $F \in R_{t-i}$.  But because $(J+\wp^{i+1})$ is
$\N^k$-homogenous, we can take $F$ to be $\N^k$-homogeneous,
and so $F = G + H$ where
$G$ is a monomial of degree $t - i$ in $x_{1,0},\ldots,x_{k,0}$
and $H \in \wp$.  Since $HM \in \wp^{i+1}$, we have 
$((J+\wp^i)/(J+\wp^{i+1}))_t = 0$ if and only if $GM \in (J + \wp^{i+1})_t$,
as desired.
\end{proof}

\begin{lemma} \label{hyperplane}
Let $P_1,\ldots,P_r,P$ be points in generic position
in $\pnk$ with $n_1 \geq \cdots \geq n_k$, and let $m_1 \geq \cdots \geq m_r$ be positive integers.
Set $J = \wp_1^{m_1} \cap \cdots \cap \wp_r^{m_r}$.  Suppose
$\underline{a} = (a_1,\ldots,a_k) \in \N^k$ is such that 
$n_k\left(\sum_{i=1}^k a_i\right) \geq \sum_{i=1}^r m_i$
and $\sum_{i=1}^k a_i \geq m_1$.  Then we can find $a_j$
hyperplanes $L_{j,1},\ldots,L_{j,a_j}$ in $\pr^{n_j}$,
that is, $ L_{j,l} \in {\bf k}[x_{j,0},\ldots,x_{j,n_j}]$ for all $l = 1, \ldots, a_j$, such
that 
\[L = \prod_{j=1}^k \left(\prod_{l=1}^{a_j} L_{j,l}\right) \in J \]
and $L$ avoids $P$.
\end{lemma}

\begin{proof}
If $r \leq n_j$ for all $j$, then for each $j$ we can find
a linear form $L_j \in {\bf k}[x_{j,0},\ldots,x_{j,n_j}]$
that passes through $P_1,\ldots,P_r$ and avoids $P$.
If we take $L_{j,l} = L_j$ for all $j$, we have
\[ L = \prod_{j=1}^k L_{j}^{a_j} \in \wp_1^{|\underline{a}|}
\cap \cdots \cap \wp_r^{|\underline{a}|} \subseteq
\wp_1^{m_1} \cap \cdots \cap \wp_r^{m_r} = J,\]
where $|\underline{a}| = \sum_{i=1}^k a_i$, since
$|\underline{a}| \geq m_1 \geq \cdots \geq m_r$.  Moreover,
$L$ avoids $P$. 

Suppose now that $n_k \leq n_{k-1} \leq \cdots
\leq n_{l+1} < r \leq n_l \leq \cdots \leq n_1$.
We shall use induction on $\sum_{i=1}^r m_i$.  Note
that if $\sum_{i=1}^r m_i \leq n_k$ then the conclusion
follows since in this case $r \leq n_k \leq n_j$ for all $j$.
If $a_k = a_{k-1} = \cdots = a_{l+1} = 0$, then the 
conclusion follows as in the case $r \leq n_j$ for all $j$.
Suppose there is $p \in \{l+1,\ldots, k\}$ such that $a_p \neq 0$.
Choose a hyperplane $L_1$ in $\pr^{n_p}$ 
($L_1 \in {\bf k}[x_{p,0},\ldots,x_{p,n_p}]$) that avoids $P$
and passes through $P_1,\ldots,P_{n_p}$.
Since $n_k(\sum_{i=1}^k a_i) \geq \sum_{i=1}^r m_i$,
we have 
\begin{eqnarray*}
n_k\left(\sum_{i=1}^k a_i\right) -n_k
 &\geq& \sum_{i=1}^r m_i - n_k  \geq \sum_{i=1}^r
m_i - n_p \\
&=& (m_1 -1) + \cdots (m_{n_p}-1) + m_{n_p+1} + \cdots + m_r.
\end{eqnarray*}
If we set $(b_1,\ldots,b_{p-1},b_p,b_{p+1},\ldots,b_k) = 
(a_1,\ldots,a_{p-1},a_p-1,a_{p+1},\ldots,a_k)$, then
we have 
\[n_k \left(\sum_{i=1}^k b_i\right)
=n_k \left(\sum_{i=1}^k a_i \right) -n_k \geq  
(m_1 -1) + \cdots (m_{n_p}-1) + m_{n_p+1} + \cdots + m_r.\]  
By induction
there exists $L_{j,1},\ldots,L_{j,b_j}$ in $\pr^{n_j}$ for
all $j$ that avoids $P$ such that 
\[L = \prod^k_{j=1}\left( \prod_{l=1}^{b_j} L_{j,l}\right)
\in \wp_1^{m_1-1}\cap \cdots \cap \wp_{n_p}^{m_{n_p}-1} \cap
\wp_{n_p+1}^{m_{n_p+1}} \cap \cdots \cap \wp_r^{m_r}.\]
If we take $L\cdot L_1$ we have the conclusion since $L_1 \in
\wp_1 \cap \cdots \cap \wp_{n_p}$ (the $a_p$ hyperplanes in $\pr^{n_p}$ are $L_{p,1}, \ldots, L_{p, b_p}$ and $L_1$).
\end{proof}

\begin{proposition} \label{riprop}
Let $P_1,\ldots,P_r,P$ be points in generic position in $\pnk$
with $n_1 \geq \cdots \geq n_k$.
Suppose $m_1 \geq \cdots \geq m_r \geq a$ are positive integers.
Set $J = \wp^{m_1} \cap \cdots \cap \wp_r^{m_r}$.  Let $t$
be the least integer such that $n_kt \geq \sum_{i=1}^r m_i + a - 1$. Then
\[ \ri(R/(J+\wp^a)) \leq \max\{m_1 + a - 1, t\}.\]
\end{proposition}

\begin{proof} Without loss of generality take $P = [1:0:\cdots:0]
\times \cdots \times [1:0:\cdots:0]$.  Then
$\wp = (x_{1,1},\ldots,x_{1,n_1},\ldots,x_{k,1},\ldots,x_{k,n_k})$.
If $r \leq n_j$ for all $j$, then we can find a hyperplane $L_j$
in $\pr^{n_j}$, i.e., $L_j \in {\bf k}[x_{j,0},\ldots,x_{j,n_j}]$, containing
$P_1,\ldots,P_r$ and avoids $P$ for each $j$.  Then $L_j 
\in \wp_1 \cap \cdots \cap \wp_r$ for all $j$.  

Suppose $G = x_{1,0}^{a_1}\cdots x_{k,0}^{a_k}$ is a monomial of degree
$m_1$ in $\{x_{1,0},\ldots,x_{k,0}\}$.  Then 
$L := L_1^{a_1}\cdots L_k^{a_k} \in \wp_1^{m_1} \cap \cdots \cap
\wp_r^{m_1} \subseteq \wp_1^{m_1} \cap \cdots \cap \wp_r^{m_r} =J$.  
We can rewrite $L_j = x_{j,0} + H_j$ where $H_j \in (x_{j,1},\ldots,x_{j,n_j})
\subseteq \wp$.  Then $L \in J$ implies $G \in J + \wp$.  Thus,
for any monomial $M$ of degree $i$ in $\wp^i$ for some $0 \leq i \leq a-1$,
$GM \in J+\wp^{i+1}$.  Since $a-1 \geq i$, this implies that
for any monomial $\overline{G}$ of degree $m_1 + a -1 - i$
in $\{x_{1,0},\ldots,x_{k,0}\}$, and any monomial $M$ of degree $i$
in $\wp^i$, $\overline{G}M \in (J + \wp^{i+1})$ because $\overline
{G}$ is divisible by a monomial of degree $m_1$.  By Lemma \ref{bounds},
this implies that $\ri(R/(J+\wp^a)) \leq m_1 + a -1$.

Suppose now that $r > n_k$.  Since $n_1 \ge \ldots \ge n_k$, by a change of coordinates we may assume that
\begin{eqnarray*}
P_1 & = & [0:1:0:\cdots:0] \times [0:1:0:\cdots:0] \times \cdots \times
[0:1:0:\cdots:0]\\
& \vdots & \\
P_{n_k} & = & [\underbrace{0:\cdots:0}_{n_k}:1:0:\cdots:0] \times [\underbrace{0:\cdots:0}_{n_k}:1:0:\cdots:0] 
\times \cdots \times [0:\cdots:0:1]
\end{eqnarray*}
So for $0 \leq j \leq n_k$, 
$\wp_{j} = (\{x_{l,q} ~|~ l=1,\ldots,k, ~q \neq j\})$. 

Let $h = \max\{m_1 + a -1, t\}$ and $0 \leq i \leq a-1$.  Suppose
now that $G = x_{1,0}^{a_1}\cdots x_{k,0}^{a_k}$ is a monomial of 
degree $h-i$ in $\{x_{1,0},\ldots,x_{k,0}\}$, and $M = \prod_{l=1}^k \prod_{q\neq 0} x_{l,q}^{c_{l,q}}$ is a monomial of degree $i$ in $\wp^i$. Because of Lemma \ref{bounds} we need to show that $GM \in (J+\wp^{i+1}).$

It can be seen that
\[M \in \wp_1^{i-\sum_{l=1}^k c_{l,1}} \cap \wp_2^{i - \sum_{l=1}^k c_{l,2}}
\cap \cdots \cap \wp_{n_k}^{i- \sum_{l=1}^k c_{l,n_k}}.\]
We also have, since $i \leq a-1,$
\[ \sum_{i=1}^k a_i = h - i \geq m_1 \geq \max\left\{ 
m_1 -  i + \sum_{i=1}^k c_{l,1},
\ldots,
m_{n_k} -  i + \sum_{i=1}^k c_{l,n_k}
\right\} \]
and 
\begin{eqnarray*}
n_k\left(\sum_{j=1}^k a_j\right) & = & n_k(h-i) = n_kh - in_k \\
& \geq & \sum_{j=1}^r m_j + a - 1 - in_k \geq \sum_{j=1}^r m_j + i  - in_k \\
& \geq & \sum_{j=1}^r m_j + \sum_{l=1}^k\sum_{q=1}^{n_k} c_{l,q} - in_k \\
& = & (m_1 - i + \sum_{l=1}^k c_{l,1}) +
\cdots + (m_{n_k} - i + \sum_{l=1}^k c_{l,n_k}) +
m_{n_k+1} + \cdots + m_r.
\end{eqnarray*}
Using Lemma \ref{hyperplane}, there exists $L_{j,1},\ldots,L_{j,a_j}
\in {\bf k}[x_{j,0},\ldots,x_{j,n_j}]$ for each $1 \leq j \leq k$ such that 
\[L = \prod_{j=1}^k \left(\prod_{q=1}^{a_j} L_{j,q}\right)
\in \wp_1^{m_1 - i + \sum c_{l,1}} \cap \cdots \cap
\wp_{n_k}^{m_{n_k} - i + \sum c_{l,n_k}} \cap \wp_{n_k+1}^{m_{n_k+1}} \cap
\cdots \cap \wp_r^{m_r},\]
and $L$ avoids $P$. This implies that $LM \in J$.

Since $L_{j,q}$ avoids $P$ we can write $L_{j,q} = x_{j,0} + H_{j,q}$
where $H_{j,q} \in (x_{j,1},\ldots,x_{j,n_j}) \subseteq \wp$.
Then $L = x_{1,0}^{a_1}\cdots x_{k,0}^{a_k} + N$ where $N \in \wp$.
Thus, since $LM \in J$, then $GM \in (J+\wp^{i+1})$ which is what we 
need to prove.
\end{proof}

\begin{theorem}  \label{mainresult2}
Suppose $P_1,\ldots,P_s$ are points in generic position
in $\pnk$ ($s \geq 2$ and $n_1 \ge \ldots \ge n_k$), and $m_1\geq m_2 \geq \cdots \geq m_s$ are positive
integers.  Set $I = \wp_1^{m_1} \cap \cdots \cap \wp_s^{m_s}$.  Then
\[ \ri(R/I) \leq \max\left\{m_1 + m_2 -1, \left[\frac{\sum_{i=1}^s m_i + n_k
-2}{n_k} \right] \right\}\]
where $[q]$ denotes the floor function.
\end{theorem}

\begin{proof}
Note that $n_1 \geq \cdots \geq n_k$, so
\[\left[ \frac{\sum_{i=1}^s m_i + n_k - 2}{n_k}\right]
= 
\max\left\{ \left[ \frac{\sum_{i=1}^s m_i + n_j - 2}{n_j}\right]\right\}_{j=1}^k
\]
Also, $\min\{t ~|~ n_kt \geq q\} = \left[ \frac{q+n_k-1}{n_k}\right]$.
So, if we take $q = \sum_{i=1}^r m_i + m_{r+1} -1$ and use
Proposition \ref{riprop} and 
induction successively, along with Lemma \ref{riformanypoints} 
we will have the conclusion.
\end{proof}

We obtain an immediate corollary which gives a bound on the regularity of the defining ideal of a scheme of fat points in $\pnk$.
\begin{corollary} \label{mainreg}
With the hypotheses as in Theorem \ref{mainresult2} we have
\[\reg(I) \leq \max\left\{m_1 + m_2 -1, \left[\frac{\sum_{i=1}^s m_i + n_k
-2}{n_k} \right] \right\} + k. \]
\end{corollary}

\begin{remark} When $k=1$ we recover the result of \cite{CTV} which was proved to be sharp. Thus, our bound in Corollary \ref{mainreg} is sharp.
\end{remark}


\end{document}